%
%
%

\magnification=1200
\pretolerance=500 \tolerance=1000 \brokenpenalty=5000
\hsize=13.5cm
\vsize=19cm
\parskip3pt plus 1pt
\parindent=0.6cm

\font\bfseventeen=cmbx10 at 17.28pt
\font\bftwelve=cmbx10 at 12pt

\font\tensans=cmss10
\font\fivesans=cmss10 at 5pt
\font\sevensans=cmss10 at 7pt
\newfam\sansfam
\textfont\sansfam=\tensans\scriptfont\sansfam=\sevensans
\scriptscriptfont\sansfam=\fivesans

\font\tenCal=eusm10
\font\sevenCal=eusm7
\font\fiveCal=eusm5
\newfam\Calfam
  \textfont\Calfam=\tenCal
  \scriptfont\Calfam=\sevenCal
  \scriptscriptfont\Calfam=\fiveCal

\font\tenBbb=msbm10
\font\sevenBbb=msbm7
\font\fiveBbb=msbm5
\newfam\Bbbfam
  \textfont\Bbbfam=\tenBbb
  \scriptfont\Bbbfam=\sevenBbb
  \scriptscriptfont\Bbbfam=\fiveBbb
\def\Bbb{\fam\Bbbfam\tenBbb}

\def\bC{{\Bbb C}}
\def\bN{{\Bbb N}}
\def\bP{{\Bbb P}}
\def\bQ{{\Bbb Q}}
\def\bR{{\Bbb R}}
\def\bZ{{\Bbb Z}}
\def\bOne{{\bf 1}}

\def\square{{\hfill \hbox{
\vrule height 1.453ex  width 0.093ex  depth 0ex
\vrule height 1.5ex  width 1.3ex  depth -1.407ex\kern-0.1ex
\vrule height 1.453ex  width 0.093ex  depth 0ex\kern-1.35ex
\vrule height 0.093ex  width 1.3ex  depth 0ex}}}
\def\qed{\phantom{$\quad$}\hfill$\square$\medskip}
\def\srelbar{\vrule width0.6ex height0.65ex depth-0.55ex}
\def\merto{\mathrel{\srelbar\kern1.3pt\srelbar\kern1.3pt\srelbar
    \kern1.3pt\srelbar\kern-1ex\raise0.28ex\hbox{${\scriptscriptstyle>}$}}}

\def\section#1|{\par\vskip .5cm\penalty -100 
\vbox{\noindent{\bftwelve #1}
\vskip 5pt}
\penalty 500}

\def\bibitem#1&#2&#3&#4&%
{\hangindent=1.66cm\hangafter=1
\noindent\rlap{\hbox{\bf #1}}\kern1.66cm{\rm #2}{\it #3}{\rm #4}} 

\long\def\claim#1|#2\endclaim{\par\vskip 5pt\noindent 
{\bf #1.}\ {\it #2}\par\vskip 5pt}

\def\today{\ifcase\month\or
January\or February\or March\or April\or May\or June\or July\or August\or
September\or October\or November\or December\fi \space\number\day,
\number\year}

\def\dbar{\overline\partial}
\def\ddbar{\partial\overline\partial}
\def\codim{\mathop{\rm codim}}
\def\Tr{\mathop{\rm Tr}\nolimits}

\def\Herm{\mathop{\rm Herm}}
\def\Vect{\mathop{\rm Vect}}
\def\mod{\mathop{\rm mod}}

\def\Alb{\mathop{\rm Alb}}

\def\PSH{\mathop{\rm PSH}}
\def\Ricci{\mathop{\rm Ricci}}

\let\ol=\overline

\def\semidirect{\mathop{\kern2pt\vrule depth-0.3pt height4.3pt 
\kern-2pt\times}\nolimits}

\def\UU{{\rm U}}
\def\SU{{\rm SU}}
\def\Sp{{\rm Sp}}
\def\rA{{\rm A}}

\def\cF{{\cal F}}
\def\cL{{\cal L}}
\def\cO{{\cal O}}
\def\cX{{\cal X}}
\def\cY{{\cal Y}}
\def\cZ{{\cal Z}}

\centerline{\bfseventeen Rationally connected manifolds and} \medskip
\centerline{\bfseventeen semipositivity of the Ricci curvature}
\bigskip\bigskip
\line{\bftwelve\hfil
Fr\'ed\'eric Campana, Jean-Pierre Demailly, Thomas Peternell\hfil} 
\vskip1cm

\noindent{\bf Abstract.} This work establishes a structure theorem for
compact K\"ahler manifolds with semipositive anticanonical bundle. 
Up to finite \'etale cover, it is proved that such manifolds split
holomorphically and isometrically as a product of Ricci flat varieties
and of rationally connected manifolds. The proof is based on a
characterization of rationally connected manifolds through the non-existence
of certain twisted contravariant tensor products of the tangent bundle, 
along with a generalized holonomy principle for pseudoeffective
line bundles. A crucial ingredient for this is the characterization
of uniruledness by the property that the anticanonical bundle is not
pseudoeffective.
\medskip

\noindent{\bf Keywords.} Compact K\"ahler manifold, anticanonical
bundle, Ricci curvature, uniruled variety, rationally connected
variety, Bochner formula, holonomy principle, fundamental group,
Albanese mapping, pseudoeffective line bundle, De Rham splitting
theorem \medskip

\noindent{\bf MSC classification 2010.} 14M22, 14J32, 32J27
\vskip30pt

\section 1. Main results|

The goal of this work is to understand the geometry of compact
K\"ahler manifolds with semipositive Ricci curvature, and especially
to study the relations that tie Ricci semipositivity with rational
connectedness. Many of the ideas are borrowed from [DPS96] and
[BDPP].  Recall that a compact complex manifold $X$ is said to be
rationally connected if any two points of $X$ can be joined by a chain
of rational curves. A line bundle $L$ is said to be
hermitian semipositive if it can be equipped with a smooth hermitian
metric of semipositive curvature form. A sufficient condition for
hermitian semipositivity is that some multiple of $L$ is spanned by
global sections; on the other hand, the hermitian semipositivity
condition implies that $L$ is numerically effective (nef) in the sense
of [DPS94], which, for $X$ projective algebraic, is equivalent to
saying that $L\cdot C\geq0$ for every curve $C$ in~$X$.  Examples
contained in [DPS94] show that all three conditions are different
(even for $X$ projective algebraic). The Ricci curvature is the curvature
of the anticanonical bundle $K_X^{-1}=\det(T_X)$, and by Yau's solution 
of the Calabi conjecture (see [Aub76], [Yau78]), a compact K\"ahler
manifold $X$ has a hermitian semipositive anticanonical bundle $K_X^{-1}$
if and only if $X$ admits a K\"ahler metric $\omega$ with $\Ricci(\omega)\geq
0$. A classical example of projective surface with $K_X^{-1}$ nef is the
complex projective plane $\bP^2_\bC$ blown-up in 9 points, no 3 of
which are collinear and no 6 of which lie on a conic; in that
case Brunella [Bru10] showed that there are configurations of the 9
points for which $K_X^{-1}$ admits a smooth (but non-real analytic) metric
with semipositive Ricci curvature; depending on some diophantine
condition introduced in [Ued82], there are also configurations for which
some multiple $K_X^{-m}$ of $K_X^{-1}$ is generated by sections and
others for which $K_X^{-1}$ is nef without any smooth metric. Finally,
let us recall that a line bundle $L\to X$ is said to be
pseudoeffective if here exists a singular hermitian metric $h$ on $L$
such that the Chern curvature current $T=i\Theta_{L,h}=-i\ddbar\log h$
is non-negative; equivalently, if $X$ is projective algebraic, this
means that the first Chern class $c_1(L)$ belongs to the closure of
the cone of effective $\bQ$-divisors.

We first give a criterion characterizing rationally connected manifolds
by the non-existence of sections in certain twisted tensor powers of
the cotangent bundle; this is only a minor variation of Theorem~5.2 
in~[Pet06], cf.\ also Remark 5.3 therein.

\claim 1.1. Criterion for rational connectedness|Let $X$ be a projective 
algebraic $n$-dimensional manifold. The following properties are equivalent.
\smallskip
\item {\rm (a)} $X$ is rationally connected.
\smallskip
\item {\rm (b)} For every invertible subsheaf 
$\cF\subset\Omega^p_X:=\cO(\Lambda^pT^*_X)$, $1\le p\le n$, 
$\cF$ is not pseudoeffective.
\smallskip
\item {\rm (c)} For every invertible subsheaf 
$\cF\subset\cO((T^*_X)^{\otimes p})$, $p\ge 1$, 
$\cF$ is not pseudoeffective.
\smallskip
\item {\rm (d)} For some $($resp.\ for any$)$ ample line bundle $A$ on $X$, 
there exists a constant $C_A>0$ such that
$$
H^0(X,(T^*_X)^{\otimes m}\otimes A^{\otimes k})=0\qquad
\hbox{for all $m,\,k\in\bN^*$ with $m\ge C_Ak$.}
$$
\endclaim

\claim 1.2. Remark|{\rm The proof follows easily from the uniruledness
criterion established in~[BDPP]$\,$: a non-singular projective variety
$X$ is uniruled if and only if $K_X$ is not pseudoeffective. 
A conjecture attributed to Mumford asserts
that the weaker assumption 
 $ ({\rm d}') \  H^0(X,(T^*_X)^{\otimes m})= 0 $
for all $m\ge 1$ should be sufficient to imply rational connectedness.
Mumford's conjecture can actually be proved by essentially
the same argument if one uses the abundance conjecture in place of the
more demanding uniruledness criterion from [BDPP] -- more specifically 
that $H^0(X,K_X^{\otimes m})=0$ for all $m\ge 1$ would imply uniruledness.}
\endclaim

\claim 1.3. Remark|{\rm By [DPS94], hypotheses 1.1~(b) and (c) make sense
on an arbitrary compact complex manifold and imply that
$H^0(X,\Omega^2_X)=0$. If $X$ is assumed to
be compact K\"ahler, then $X$ is automatically
projective algebraic by Kodaira [Kod54], therefore, 1.1~(b) or (c) also 
characterize rationally connected manifolds among all compact
K\"ahler ones.\qed}
\endclaim

The following structure theorem generalizes the
Bogomolov-Kobayashi-Beau\-ville structure theorem for Ricci-flat manifolds
([Bog74a], [Bog74b], [Kob81], [Bea83]) to the Ricci semi\-positive case. Recall
that a~{\it holomorphic symplectic manifold} $X$ is a compact K\"ahler 
manifold admitting a holomorphic symplectic $2$-form $\omega$ (of maximal 
rank everywhere); in particular $K_X=\cO_X$. A {\it Calabi-Yau} manifold is
a simply connected projective manifold with $K_X=\cO_X$ and
$H^0(X,\Omega^p_X)=0$ for $0<p<n=\dim X$ (or a finite \'etale quotient
of such a manifold).

\claim 1.4. Structure theorem|Let $X$ be a compact K\"ahler manifold with
$K_X^{-1}$ hermitian semipositive. Then
\smallskip
\item {\rm (a)} The universal cover $\widetilde X$ admits a holomorphic and
isometric splitting
$$\widetilde X\simeq\bC^q\times\prod Y_j\times\prod S_k\times\prod Z_\ell$$
where $Y_j$, $S_k$, and $Z_\ell$ are compact simply connected K\"ahler
manifolds of respective dimensions $n_j$, $n'_k$, $n''_\ell$ 
with irreducible holonomy, $Y_j$ being Calabi-Yau manifolds
$($holonomy $\SU(n_j))$, $S_k$ holomorphic symplectic manifolds
$($holonomy $\Sp(n'_k/2))$, and $Z_\ell$ rationally connected manifolds
with $K_{Z_\ell}^{-1}$ semipositive $($holonomy $\UU(n''_\ell)$, unless
$Z_\ell$ is hermitian symmetric of compact type$)$.
\smallskip
\item {\rm (b)} There exists a finite \'etale Galois cover $\widehat X\to X$
such that the Albanese va\-riety $\Alb(\widehat X)$ is a $q$-dimensional torus
and the Albanese map \hbox{$\alpha:\widehat X\to\Alb(\widehat X)$} is an
$($isometrically$)$ locally trivial holomorphic fiber bundle whose fibers 
are products $\prod Y_j\times\prod S_k\times\prod Z_\ell$ of the type 
described in {\rm a)}. Even more holds after possibly another finite \'etale cover\/$:$  
$\hat X$ is a fiber bundle with fiber $\prod Z_\ell$ on $\prod Y_j \times \prod S_k \times \Alb(\widehat X). $
\smallskip
\item {\rm (c)} We have $\pi_1(\widehat X)\simeq\bZ^{2q}$ and $\pi_1(X)$ is an
extension of a finite group $\Gamma$ by the normal subgroup $\pi_1(\widehat X)$.
In particular there is an exact sequence
$$0\to \bZ^{2q} \to \pi_1(X) \to \Gamma \to 0,$$
and the fundamental group $\pi_1(X)$ is almost abelian.\vskip0pt
\endclaim

The proof relies on the holonomy principle, and on De Rham's splitting theorem 
[DR52] and Berger's classification [Ber55].
Foundational background can be found in papers by Lichnerowicz [Lic67], 
[Lic71], and Cheeger-Gromoll [CG71], [CG72]. The restricted holonomy group
of a hermitian vector bundle $(E,h)$ of rank $r$ is by definition the subgroup 
$H\subset\UU(r)\simeq U(E_{z_0})$ generated by parallel transport 
operators with respect to the Chern connection $\nabla$ of $(E,h)$, 
along loops based at $z_0$ that are contractible (up to conjugation,
$H$ does not depend on the base point~$z_0$). We need here a 
generalized ``pseudoeffective'' version of the holonomy principle, which 
can be stated as follows.

\claim 1.5. Generalized holonomy principle|Let $E$ be a holomorphic
vector bundle of rank $r$ over a compact complex manifold $X$. Assume
that $E$ is equipped with a smooth hermitian structure $h$ and $X$ with
a hermitian metric~$\omega$, viewed as a smooth positive $(1,1)$-form
$\omega=i\sum \omega_{jk}(z)dz_j\wedge d\overline z_k$.
Finally, suppose that the $\omega$-trace of the Chern curvature tensor
$\Theta_{E,h}$ is semipositive, that is
$$
i\Theta_{E,h}\wedge{\omega^{n-1}\over (n-1)!}=
B\,{\omega^n\over n!},\qquad B\in\Herm(E,E),\quad \hbox{with}~~
B\ge 0~~\hbox{on $X$},
$$
and denote by $H$ the restricted holonomy group of $(E,h)$.\vskip-\parskip
\penalty-10000\smallskip

\item{\rm (a)} If there exists an invertible sheaf $\cL\subset
\cO((E^*)^{\otimes m})$ which is pseudoeffective as a line bundle, 
then $\cL$ is flat and $\cL$ is invariant under parallel transport by
the connection of $(E^*)^{\otimes m}$ induced by the Chern connection
$\nabla$ of~$(E,h)\,;$ in fact, $H$ acts trivially on~$\cL$.
\smallskip
\item{\rm (b)} If $H$ satisfies $H=\UU(r)$, then none of the 
invertible sheaves $\cL\subset\cO((E^*)^{\otimes m})$ can be pseudoeffective
for $m\ge 1$.
\vskip0pt
\endclaim

The generalized holonomy principle is based on an extension of the
Bochner formula as found in [BY53], [Ko83]$\,$: for $(X,\omega)$
K\"ahler, every section $u$ in $H^0(X,(T^*_X)^{\otimes m})$ satisfies
$$
\Delta (\Vert u \Vert^2) = \Vert\nabla u\Vert^2 + Q(u),
\leqno(1.6)
$$
where $Q(u)\geq m\lambda_1\Vert u\Vert^2$ is bounded from below by
the smallest eigenvalue  $\lambda_1$ of the Ricci curvature 
tensor of~$\omega$. If $\lambda_1\ge 0$, the equality
$\int_X\Delta(\Vert u\Vert^2)\omega^n=0$ implies $\nabla u=0$ and
$Q(u)=0$. The generalized principle consists essentially in considering
a general vector bundle $E$ rather than $E=T^*_X$, and replacing
$\Vert u\Vert^2_\omega$ with $\Vert u\Vert^2_\omega e^\varphi$
where $u$ is a local trivializing section of $\cL$, where $\varphi$ is 
the corresponding local plurisubharmonic weight representing the metric
of $\cL$ and $\omega$ a Gauduchon metric, cf.\ (3.2).

\claim 1.7. Remark|{\rm If one makes the weaker assumption that $K_X^{-1}$
is nef, then Qi Zhang [Zha96, Zha05] proved that the Albanese mapping
$\alpha:X\to\Alb(X)$ is surjective in the case where $X$ is projective,
and P\u{a}un [Pau12] recently extended this result to the general
K\"ahler case (cf.\ also [CPZ03]). One may wonder whether there still
exists a holomorphic splitting
$$
\widetilde X\simeq\bC^q\times\prod Y_j\times\prod S_k\times\prod Z_\ell
$$
of the universal covering as above. However the example where $X=\bP(E)$ 
is the ruled surface over an elliptic curve $C=\bC/(\bZ+\bZ\tau)$ associated 
with a non-trivial rank~$2$ bundle $E\to C$ with
$$
0\to\cO_C\to E\to\cO_C\to 0
$$
shows that $\widetilde X=\bC\times\bP^1$ cannot be an isometric
product for a K\"ahler metric $\omega$ on~$X$. Actually, such a
situation would imply that $K_X^{-1}=\cO_{\bP(E)}(1)$ is semipositive,
but we know by [DPS94] that $\cO_{\bP(E)}(1)$ is nef and non-semipositive. 
Under the mere assumption that $K_X^{-1}$ is nef, it is
unknown whether the Albanese map $\alpha:X\to\Alb(X)$ is a submersion,
unless $X$ is a projective threefold [PS98], 
and even if it is supposed to be so, it seems to be unknown whether the
fibers of $\alpha$ may exhibit non-trivial variation of the complex
structure (and whether they are actually products of Ricci flat manifolds
by rationally connected manifolds). The main difficulty is that, a
priori, the holonomy argument used here breaks down -- a possibility
would be to consider some sort of ``asymptotic holonomy'' for a
sequence of K\"ahler metrics satisfying $\Ricci(\omega_\varepsilon)\ge
-\varepsilon\omega_\varepsilon$, and dealing with the Gromov-Hausdorff
limit of the variety.\qed}
\endclaim

This work was completed while the three authors were visiting the
Mathematisches Forschungsinstitut Oberwolfach in September 2012.
They wish to thank the Institute for its hospitality and the 
exceptional quality of the environment.

\section 2. Proof of the criterion for rational connectedness|

In this section we prove Criterion 1.1.
Observe first that if $X$ is rationally connected, then there exists 
an immersion $f:\bP^1\subset X$ passing through any given finite subset of $X$
such that $f^*T_X$ is ample, see e.g.\ [Kol96, Theorem~3.9, p.~203]. 
In other words $f^*T_X=\bigoplus\cO_{\bP^1}(a_j)$, $a_j>0$, while
$f^*A=\cO_{\bP^1}(b)$, $b>0$. Hence
$$H^0(\bP^1,f^*((T^*_X)^{\otimes m}\otimes A^{\otimes k}))=0\qquad
\hbox{for $m>kb/\min(a_i)$.}
$$
As the immersion $f$ moves freely in $X$, we immediately see from this that
1.1~(a) implies 1.1~(d) with any constant value $C_A>b/\min(a_j)$.

To see that 1.1~(d) implies 1.1~(c), assume that $\cF \subset 
(T^*_X)^{\otimes p}$ is a  pseudoeffective line bundle. Then there exists
$k_0\gg 1$ such that
$$ H^0(X, \cF^{\otimes m}\otimes A^{k_0}) \ne 0 $$
for all $m\ge 0$ (for this, it is sufficient to take $k_0$ such that
$A^{k_0}\otimes(K_X\otimes G^{n+1})^{-1}>0$ for some very ample line 
bundle~$G$). This implies
$H^0(X,(T^*_X)^{\otimes mp} \otimes A^{k_0}) \ne 0$ for all~$m$,
contradicting assumption 1.1~(d).

The implication 1.1~(c)$~\Rightarrow~$1.1~(b) is trivial.

It remains to show that 1.1~(b) implies 1.1~(a). First note that $K_X$ 
is not pseudoeffective, as one sees by applying the assumption 1.1~(b) 
with $p=n$. Hence $X$ is uniruled by [BDPP].
We consider the quotient with maximal rationally connected
fibers (rational quotient or MRC fibration, see [Cam92], [KMM92])
$$ f: X \merto W$$
to a smooth projective variety $W$. By [GHS01], $W$ is not uniruled, otherwise we could lift the ruling to $X$ and the fibers of $f$ would not be maximal.
We may further assume that $f$ is holomorphic. In fact, assumption 1.1~(b) is 
invariant under blow-ups. To see this, let 
$\pi: \hat X \to X$ be a birational morphisms from a projective manifold  $\hat X$ and consider a line bundle $\hat \cF \subset \Omega^p_{\hat X}.$ 
Then $\pi_*(\hat \cF) \subset \pi_*(\Omega_{\hat X}^p) = \Omega^p_X, $ hence we introduce the line bundle
$$ \cF := (\pi_*(\hat \cF))^{**} \subset \Omega^p_X.$$
Now, if $\hat \cF$ were pseudoeffective, so would be $\cF$.
Thus 1.1~(b) is invariant under $\pi$ and we may suppose $f$ holomorphic. 
In order to show that $X$ is rationally connected, we need to prove that $p:=
\dim W = 0$. Otherwise $K_W = \Omega^p_W$ is pseudoeffective by [BDPP], and
we obtain a pseudo-effective invertible subsheaf
$\cF := f^*(\Omega^p_W) \subset \Omega^p_X$, in contradiction
with 1.1~(b).\qed

\section 3. Bochner formula and generalized holonomy principle|

Let $(E,h)$ be a hermitian holomorphic vector bundle over a
$n$-dimensional compact complex manifold~$X$. The semipositivity
hypothesis on $B=\Tr_\omega\Theta_{E,h}$ is invariant by a conformal
change of metric $\omega$. Without loss of generality we can assume
that $\omega$ is a Gauduchon metric, i.e.\ that
$\ddbar\omega^{n-1}=0$, cf.\ [Gau77]. We consider the Chern connection
$\nabla$ on $(E,h)$ and the corresponding parallel transport
operators. At every point $z_0\in X$, there exists a local coordinate
system $(z_1,\ldots,z_n)$ centered at $z_0$ (i.e.\ $z_0=0$ in
coordinates), and a holomorphic frame $(e_\lambda(z))_{1\le\lambda\le r}$
such that
$$
\leqalignno{
\langle e_\lambda(z),e_\mu(z)\rangle_h&=\delta_{\lambda\mu}-
\sum_{1\le j,k\le n} c_{jk\lambda\mu}z_j\ol z_k+O(|z|^3),\qquad
1\le\lambda,\mu\le r,&(3.1)\cr
\Theta_{E,h}(z_0)&=\sum_{1\le j,k,\lambda,\mu\le n} c_{jk\lambda\mu}
dz_j\wedge d\ol z_k\otimes e_\lambda^*\otimes e_\mu,\qquad
c_{kj\mu\lambda}=\ol{c_{jk\lambda\mu}},&(3.1')\cr}
$$
where $\delta_{\lambda\mu}$ is the Kronecker symbol and $\Theta_{E,h}(z_0)$
is the curvature tensor of the Chern connection $\nabla$ of $(E,h)$
at~$z_0$.

Assume that we have an invertible sheaf $\cL\subset \cO((E^*)^{\otimes m})$
that is pseudoeffective. There exist a covering $U_j$ by coordinate balls
and holomorphic sections $f_j$ of $\cL_{|U_j}$ generating $\cL$ over $U_j$.
Then $\cL$ is associated with the \v{C}ech cocycle $g_{jk}$ in~$\cO_X^*$
such that $f_k=g_{jk}f_j$, and the singular hermitian metric $e^{-\varphi}$ 
of $\cL$ is
defined by a collection of plurisubharmonic functions $\varphi_j\in\PSH(U_j)$
such that $e^{-\varphi_k}=|g_{jk}|^2e^{-\varphi_j}$. It follows that we have
a globally defined bounded measurable function 
$$
\psi=e^{\varphi_j}\Vert f_j\Vert^2=e^{\varphi_j}\Vert f_j\Vert^2_{h^{*m}}
$$
over $X$, which can be viewed also as the hermitian metric ratio 
$(h^*)^m/e^{-\varphi}$ along~$\cL$, i.e.\ $\psi=(h^*)^m_{|\cL}e^\varphi$.
We are going to compute the Laplacian $\Delta_\omega\psi$.
For simplicity of notation, we omit the index $j$ and consider a
local holomorphic section $f$ of $\cL$ and a local weight 
$\varphi\in\PSH(U)$ on some open subset $U$ of $X$. In a neighborhood of
an arbitrary point $z_0\in U$, we write
$$
f=\sum_{\alpha\in\bN^m}f_\alpha\,e^*_{\alpha_1}\otimes\ldots\otimes e^*_{\alpha_m},
\qquad f_\alpha\in\cO(U),
$$
where $(e^*_\lambda)$ is the dual holomorphic frame of $(e_\lambda)$ in 
$\cO(E^*)$. The hermitian matrix of $(E^*,h^*)$ is the
transpose of the inverse of the hermitian matrix of $(E,h)$, hence
(3.1) implies
$$
\langle e^*_\lambda(z),e^*_\mu(z)\rangle_h=\delta_{\lambda\mu}+
\sum_{1\le j,k\le n} c_{jk\mu\lambda}z_j\ol z_k+O(|z|^3),\qquad
1\le\lambda,\mu\le r.
$$
On the open set $U$ the function $\psi=(h^*)^m_{|\cL}e^\varphi$ is given by
$$
\psi=\Big(\sum_{\alpha\in\bN^m}|f_\alpha|^2+
\sum_{\alpha,\beta\in\bN^m,\,1\le j,k\le n,\,1\le\ell\le m}
f_\alpha\overline{f_\beta}\,
c_{jk\beta_\ell\alpha_\ell}z_j\overline z_k+O(|z|^3)|f|^2\Big)
e^{\varphi(z)}.
$$
By taking $i\ddbar(...)$ of this at $z=z_0$ in the sense of distributions 
(that is, for almost every $z_0\in X$), we find
$$
\eqalign{
i\ddbar\psi=e^\varphi\Big(|f|^2
i\ddbar\varphi&+i\langle
\partial f+f\partial\varphi,\partial f+f\partial\varphi\rangle
+{}\cr
&+\sum_{\alpha,\beta,j,k,1\le\ell\le m}f_\alpha\overline{f_\beta}\,
c_{jk\beta_\ell\alpha_\ell}\,idz_j\wedge d\overline z_k\Big).
\cr}
$$
Since $i\ddbar\psi\wedge{\omega^{n-1}\over (n-1)!}=
\Delta_\omega\psi\,{\omega^n\over n!}$ (we actually take this as
a definition of $\Delta_\omega$), a multiplication by $\omega^{n-1}$ 
yields the fundamental inequality
$$
\Delta_\omega\psi\ge |f|^2e^\varphi(\Delta_\omega\varphi+m\lambda_1)+
|\nabla^{1,0}_hf+f\partial\varphi|^2_{\omega,h^{*m}}~e^\varphi\leqno(3.2)
$$
where $\lambda_1(z)\ge 0$ is the lowest eigenvalue of the hermitian
endomorphism \hbox{$B=\Tr_\omega\Theta_{E,h}$} at an arbitrary point $z\in X$.
As $\ddbar\omega^{n-1}=0$, we have
$$
\int_X\Delta\psi{\omega^n\over n!}=
\int_Xi\ddbar\psi\wedge{\omega^{n-1}\over (n-1)!}=
\int_X\psi\wedge{i\ddbar(\omega^{n-1})\over (n-1)!}=0
$$
by Stokes' formula. Since $i\ddbar\varphi\ge 0$, (3.2) implies
$\Delta_\omega\varphi=0$, i.e.\ $i\ddbar\varphi=0$, and
$\nabla^{1,0}_hf+f\partial\varphi=0$ almost everywhere. This means in
particular that the line bundle $(\cL,e^{-\varphi})$ is flat.
In each coordinate ball $U_j$ the pluriharmonic function $\varphi_j$
can be written $\varphi_j=w_j+\overline w_j$ for some holomorphic
function $w_j\in\cO(U_j)$, hence $\partial\varphi_j=dw_j$ and the condition 
$\nabla^{1,0}_hf_j+f_j\partial\varphi_j=0$ can be rewritten
$\nabla^{1,0}_h(e^{w_j}f_j)=0$ where $e^{w_j}f_j$ is a local holomorphic 
section. This shows that $\cL$ must be invariant by parallel transport
and that the local holonomy of the Chern connection of $(E,h)$ acts 
trivially on $\cL$. Statement 1.5~(a) follows.

Finally, if we assume that the restricted holonomy group $H$ of
$(E,h)$ is equal to~$\UU(r)$, there cannot exist any holonomy
invariant invertible subsheaf $\cL\subset\cO((E^*)^{\otimes m})$,
$m\ge 1$, on which $H$ acts trivially, since the natural
representation of $\UU(r)$ on $(\bC^r)^{\otimes m}$ has no invariant
line on which $\UU(r)$ induces a trivial action.  Property~1.5~(b) is
proved.\qed

\section 4. Proof of the structure theorem| 

We suppose here that $X$ is equipped with a K\"ahler metric $\omega$ such that
${\rm Ricci}(\omega)\geq 0$, and we set $n=\dim_{\bC}X$. We consider the
holonomy representation of the tangent bundle $E=T_X$ equipped with
the hermitian metric $h=\omega$. Here
$$
B=\Tr_\omega\Theta_{E,h}=\Tr_\omega\Theta_{T_X,\omega}\ge 0
$$
is nothing but the Ricci operator.
\medskip

\noindent {\it Proof of\/}~1.4~(a). Let $$(\widetilde
X,\omega)\simeq\prod(X_i,\omega_i)$$ be the De Rham decomposition of
$(\widetilde X,\omega)$, induced by a decomposition of the holonomy
representation in irreducible representations. Since the holonomy is
contained in $\UU(n)$, all factors $(X_i,\omega_i)$ are K\"ahler
manifolds with irreducible holonomy and holonomy group $H_i\subset
\UU(n_i)$, $n_i=\dim X_i$. By Cheeger-Gromoll [CG71], there is
possibly a flat factor $X_0=\bC^q$ and the other factors $X_i$, $i\ge
1$, are compact and simply connected. Also, the product structure
shows that each $K_{X_i}^{-1}$ is hermitian semi\-positive.  By
Berger's classification of holonomy groups [Ber55] there are only
three possibilities, namely $H_i=\UU(n_i)$, $H_i=\SU(n_i)$ or
$H_i=\Sp(n_i/2)$, unless $X_i$ is a Hermitian symmetric space, then
necessarily of compact type; such symmetric spaces have been classified
by E.~Cartan, they are rational homogeneous, hence rationally connected
(their holonomy groups are also well known, see e.g.\ [Bes87, \S$\,$10]).
The case $H_i=\SU(n_i)$ leads to $X_i$ being a
Calabi-Yau manifold, and the case $H_i=\Sp(n_i/2)$ implies that $X_i$
is holomorphic symplectic (see e.g.\ [Bea83]). Now, if $H_i=\UU(n_i)$,
the generalized holonomy principle~1.5 shows that none of the
invertible subsheaves $\cL\subset\cO((T^*_{X_i})^{\otimes m})$ can be
pseudoeffective for $m\ge 1$.  Therefore $X_i$ is rationally connected
by Criterion~1.1.\medskip

\noindent {\it Proof of\/}~1.4~(b). Set $X'=\smash{\prod_{i\geq 1}X_i}$.
The group of covering transformations acts on the product
$\smash{\widetilde X}=\bC^q\times X'$ by holomorphic isometries of the
form $x=(z,x')\mapsto(u(z),v(x'))$. At this point, the argument is
slightly more involved than in Beauville's paper [Bea83], because the
group $G'$ of holomorphic isometries of $X'$ need not be finite ($X'$
may be for instance a projective space); instead, we imitate the proof
of ([CG72], Theorem 9.2) and use the fact that $X'$ and $G'={\rm
  Isom}(X')$ are compact. Let $E_q=\bC^q\semidirect U(q)$ be the group
of unitary motions of $\bC^q$. Then $\pi_1(X)$ can be seen as a
discrete subgroup of $E_q\times G'$. As $G'$ is compact, the kernel of
the projection map $\pi_1(X)\to E_q$ is finite and the image of
$\pi_1(X)$ in $E_q$ is still discrete with compact quotient. This
shows that there is a subgroup $\Gamma$ of finite index in $\pi_1(X)$
which is isomorphic to a crystallographic subgroup of $\bC^q$. By
Bieberbach's theorem, the subgroup $\Gamma_0\subset\Gamma$ of elements
which are translations is a subgroup of finite index. Taking the
intersection of all conjugates of $\Gamma_0$ in $\pi_1(X)$, we find a
normal subgroup $\Gamma_1\subset\pi_1(X)$ of finite index, acting by
translations on $\bC^q$. Then $\widehat X=\widetilde X/\Gamma_1$ is a
fiber bundle over the torus $\bC^q/\Gamma_1$ with $X'$ as fiber and
$\pi_1(X')=1$. Therefore $\widehat X$ is the desired finite \'etale
covering of~$X$.  \smallskip
For the second assertion we consider fiberwise the rational quotient 
and obtain a factorization 
$$ \widehat X \mathrel{\buildrel {\beta}  \over {\to }} W 
\mathrel{\buildrel {\gamma} \over {\to }} \Alb(\widehat X) $$
with fiber bundles $\beta $ (fiber $\prod Z_\ell)$ and $\gamma$ 
(fiber $\prod Y_j \times \prod S_k$). Since clearly $K_W \equiv 0,$
the claim follows from the Beauville-Bogomolov decomposition theorem.
\medskip

\noindent {\it Proof of\/}~1.4~(c). The statement is an immediate
consequence of 1.4~(b), using the homo\-topy exact sequence of a fibration.\qed

\section 5. Further remarks| 

We finally point out two direct consequences of Theorem 1.4. 
Since the property 
$$ H^0(X,(T_X^*)^{\otimes m}) = 0\qquad (m \geq 1)$$
is invariant under finite \'etale covers, 
we obtain immediately from Theorem  1.4$\,:$

\claim 5.1. Corollary| Let $X$ be a compact K\"ahler manifold with 
$K_X^{-1}$ hermitian semi-positive. Assume that $H^0(X,(T_X^*)^{\otimes m}) = 0$
for all positive $m$. Then $X$ is rationally connected.
\endclaim

\noindent This establishes Mumford's conjecture in case $X$ has
semi-positive Ricci curvature.

\noindent Theorem 1.4 also gives strong implications for small
deformations of a manifold with semi-positive Ricci curvature:

\claim 5.2. Corollary| Let $X$ be a compact K\"ahler manifold with
$K_X^{-1}$ hermitian semi-positive. Let $\pi: \cX \to \Delta $ be
a proper submersion from a K\"ahler manifold $\cX$ to the unit disk
$\Delta \subset \bC$.  Assume that $X_0 = \pi^{-1}(0) \simeq X$. Then
there exists a finite \'etale cover $\widehat \cX \to \cX$
with projection $\hat \pi: \widehat \cX \to \Delta$ such that -
after possibly shrinking $\Delta$ - the following holds.  \smallskip

\item{\rm(a)} The relative Albanese map $\alpha: \widehat \cX \to
{\rm Alb}(\cX / \Delta) $ is a surjective submersion; thus the
Albanese map $\alpha_t: \widehat X_t = \widehat \pi^{-1}(t) \to {\rm Alb}(X_t)$ is a surjective submersion for all $t$.

\item{\rm(b)} Every fiber of $\alpha_t$ is a product of Calabi-Yau
manifolds, irreducible symplectic manifolds and irreducible
rationally connected manifolds.
\item{\rm(c)} There exists a factorization of $\alpha\,:$
$$  \widehat \cX \mathrel{\buildrel{\beta} \over {\to}}\cY 
\mathrel{\buildrel{\gamma} \over {\to}} {\rm Alb}(\cX/\Delta)   $$
such that $\beta_t = \beta_{|\smash{\widehat X_t}}$ is a submersion and a
rational quotient of $\widehat X_t$ for all $t$, and $\gamma_t =
\gamma_{|Y_t}$ is a trivial fiber bundle.\vskip0pt
\endclaim 

Corollary 5.2 is an immediate consequence of Theorem 1.4 and the following proposition.

\claim 5.3. Proposition| Let $\pi:\cY \to \Delta $ be a proper
K\"ahler submersion over the unit disk. Assume that $Y_0 \simeq \prod
X_i \times \prod Y_j \times \prod Z_k$ with $X_i$ Calabi-Yau, $Y_j$
irreducible symplectic and $Z_k$ irreducible rationally
connected. Then $($possibly after shrinking $\Delta)$ every $Y_t$ has
a decomposition
$$  Y_t \simeq \prod X_{i,t} \times \prod Y_{j,t} \times \prod Z_{k,t}$$
with factors of the same types as above, and the factors form families
$\cX_i,$ $\cY_j$ and $\cZ_k$.
\endclaim 

\noindent {\it Proof.} It suffices to treat the case of two factors,
say $Y_0 = A_1 \times A_2$ where the $A_i$ are Calabi-Yau, irreducible
symplectic or rationally connected.  Since $H^1(A_j,\cO_{A_j})=0$,
the factors $A_j$ deform to the neighboring $Y_t.$ By the
properness of the relative cycle space we obtain families $q_i: U_i
\to S_i$ over $\Delta$ with projections $p_i: U_i \to \cY$.
Possibly after shrinking~$\Delta$, this yields holomorphic maps $f_i:
\cY \to S_i$. Then the map
$$ f_1 \times f_2: \cY \to S_1 \times S_2 $$
is an isomorphism, since $A_t \cdot B_t = A_0 \cdot B_0 = 1$. This
gives the families $(A_i)_t$ we are looking for.\qed
\vfill\eject

\section Appendix: a flag variety version of the holonomy principle|

\centerline{by Jean-Pierre Demailly, Institut Fourier}
\medskip

Our goal here is to derive a related version of the holonomy principle
over flag varieties, based on a modified Bochner formula which we hope
to be useful in other contexts (especially since no assumption on the
base manifold is~needed). If $E$ is as before a holomorphic vector
bundle of rank $r$ over a $n$-dimensional complex manifold, we denote
by $F(E)$ the flag manifold of~$E$, namely the bundle $F(E)\to X$
whose fibers consist of flags
$$
\xi:~~E_x=V_0\supset V_1\supset \ldots \supset V_r=\{0\},\quad\dim E_x=r,
\quad\codim V_\lambda=\lambda,
$$ 
in the fibers of~$E$, along with the natural projection 
$\pi:F(E)\to X$, $(x,\xi)\mapsto x$.
We let~$Q_\lambda$, $1\le \lambda\le r$ be the tautological line bundles over
$F(E)$ such that
$$
Q_{\lambda,\xi}=V_{\lambda-1}/V_\lambda,
$$
and for a weight $a=(a_1,\ldots,a_r)\in\bZ^r$ we set
$$
Q^a=Q_1^{a_1}\otimes\ldots\otimes Q_r^{a_r}.
$$
In additive notation, viewing the $Q_j$ as divisors, we also denote
$$
a_1Q_1+\ldots+a_rQ_r
$$
any real linear combination ($a_j\in\bR$). Our goal is to compute
explicitly the curvature tensor of the line bundles $Q^a$ with respect
to the tautological metric induced by~$h$. For convenience of
notation, we prefer to work on the dual flag manifold~$F(E^*)$,
although there is a biholomorphism $F(E)\simeq F(E^*)$ given by
$$
(E_x=W_0\supset W_1\supset \ldots \supset W_r=\{0\})~~\longmapsto~~
(E^*_x=V_0\supset V_1\supset \ldots \supset V_r=\{0\})
$$
where $V_\lambda=W_{r-\lambda}^\dagger$ is the orthogonal subspace of
$W_{r-\lambda}$ in~$E^*_x$. In this context, we have an isomorphism
$$
V_{\lambda-1}/V_\lambda=W_{r-\lambda+1}^\dagger/W_{r-\lambda}^\dagger\simeq
(W_{r-\lambda}/W_{r-\lambda+1})^*.
$$
This shows that $Q^a\to F(E^*)$ is isomorphic to $Q^b\to F(E)$ where
$b_\lambda=-a_{r-\lambda+1}$, that is
$$
(b_1,b_2,\ldots,b_{r-1},b_r)=(-a_r,-a_{r-1},\ldots,-a_2,-a_1).
$$
We now proceed to compute the curvature of $Q^a\to F(E^*)$, 
using the same notation as in section~3. In a neighborhood of every point
$z_0\in X$, we can find a local coordinate system $(z_1,\ldots,z_n)$
centered at $z_0$ and a holomorphic frame $(e_\lambda)_{1\le\lambda\le r}$
such that
$$
\leqalignno{
\langle e_\lambda(z),e_\mu(z)\rangle&=\bOne_{\{\lambda=\mu\}}-
\sum_{1\le j,k\le n} c_{jk\lambda\mu}z_j\ol z_k+O(|z|^3),\qquad
1\le\lambda,\mu\le r,&(\rA.1)\cr
\Theta_{E,h}(z_0)&=\sum_{1\le j,k,\lambda,\mu\le n} c_{jk\lambda\mu}
dz_j\wedge d\ol z_k\otimes e_\lambda^*\otimes e_\mu,\qquad
c_{kj\mu\lambda}=\ol{c_{jk\lambda\mu}},&(\rA.1')\cr}
$$
where $\bOne_S$ denotes the characteristic function of the set~$S$.
For a given point $\xi_0\in F(E^*_{z_0})$ in the flag
variety, one can always adjust the frame $(e_\lambda)$ in such a way
that the flag corresponding to $\xi_0$ is given by
$$
V_{\lambda,0}=\Vect(e_1,\ldots,e_\lambda)^\dagger\subset E^*_{z_0}.
\leqno(\rA.2)
$$
A point $(z,\xi)$ in a neighborhood of $(z_0,\xi_0)$ is likewise
represented by the flag associated with the holomorphic tangent 
frame $(\widetilde e_\lambda(z,\xi))_{1\le\lambda\le r}$ defined by
$$
\widetilde e_\lambda(z,\xi)=e_\lambda(z)+
\sum_{\lambda<\mu\le r}\xi_{\lambda\mu}e_\mu(z),\qquad
(\xi_{\lambda\mu})_{1\le\lambda<\mu\le r}\in\bC^{r(r-1)/2}.\leqno(\rA.3)
$$
We obtain in this way a local coordinate system $(z_j,\xi_{\lambda\mu})$ near
$(z_0,\xi_0)$ on the total space of $F(E^*)$, where the 
$(\xi_{\lambda\mu})$ are the fiber coordinates. The frame
$\widetilde e(z,\xi)$ is not orthonormal, but by the 
Gram-Schmidt orthogonalization process, the flag $\xi$ is also induced
by the (non-holomorphic) orthonormal frame $(\widehat e_\lambda(z,\xi))$
obtained inductively by putting $\widehat e_1=\widetilde e_1/|\widetilde e_1|$
and
$$
\widehat e_\lambda=\Big(\widetilde e_\lambda-\sum_{1\le\mu<\lambda}
\langle\widetilde e_\lambda,\widehat e_\mu\rangle\,\widehat e_\mu\Big)/
\hbox{(norm of numerator)}.
$$
Straightforward calculations imply that the hermitian inner products
involved are $O(|\xi|+|z|^2)$ and the norms equal to $1+O((|\xi|+|z|)^2)$, 
hence we get
$$
\widehat e_\lambda(z,\xi)=e_\lambda(z,\xi)+
\sum_{\lambda<\mu\le r}\xi_{\lambda\mu}e_\mu(z)
-\sum_{1\le\mu<\lambda}\ol\xi_{\mu\lambda}e_\mu(z)+O\big((|\xi|+|z|)^2\big),
$$
and more precisely (omitting variables for simplicity of notation)
$$
\leqalignno{
\widehat e_\lambda
&=\Big(1-{1\over 2}\sum_{1\le\mu<\lambda}|\xi_{\mu\lambda}|^2
-{1\over 2}\sum_{\lambda<\mu\le r}|\xi_{\lambda\mu}|^2
+{1\over 2}\sum_{1\le j,k\le n}
c_{jk\lambda\lambda}z_j\ol z_k\Big)e_\lambda
+\sum_{\lambda<\mu\le r}\xi_{\lambda\mu}e_\mu\cr
&\kern35pt{}-\sum_{1\le\mu<\lambda}\Big(\ol\xi_{\mu\lambda}
+\sum_{\lambda<\nu\le r}\xi_{\lambda\nu}\ol\xi_{\mu\nu}
-\sum_{1\le j,k\le n}c_{jk\lambda\mu}z_j\ol z_k\Big)e_\mu
+O\big((|\xi|+|z|)^3\big).&(\rA.4)\cr}
$$
The curvature of the tautological line bundle $Q_\lambda=V_{\lambda-1}/V_\lambda$
can be evaluated by observing that the dual line bundle
$$
Q_\lambda^*=V_\lambda^\dagger/V_{\lambda-1}^\dagger=
\Vect(\widetilde e_1,\ldots,\widetilde e_\lambda)/
\Vect(\widetilde e_1,\ldots,\widetilde e_{\lambda-1})
$$
admits a holomorphic section given by
$$
v_\lambda(z,\xi)=\widetilde e_\lambda(z,\xi)
\mod\Vect(\widetilde e_1,\ldots,\widetilde e_{\lambda-1}).
$$
The tautological norm of this section is
$$
\eqalign{
|v_\lambda|^2&=|\widetilde e_\lambda|^2-\sum_{1\le\mu<\lambda}
|\langle\widetilde e_\lambda,\widehat e_\mu\rangle|^2\cr
&=1-\sum_{1\le j,k\le n}c_{jk\lambda\lambda}z_j\ol z_k+
\sum_{\lambda<\mu\le r}
|\xi_{\lambda\mu}|^2-\sum_{1\le\mu<\lambda}|\xi_{\mu\lambda}|^2
+O\big((|z|+|\xi|)^3\big).\cr}
$$
Therefore we obtain the formula
$$
\eqalign{
\Theta_{Q_\lambda}(z_0,\xi_0)&=\ddbar\log|v_\lambda|^2_{|(z_0,\xi_0)}\cr
&=-\sum_{1\le j,k\le n}c_{jk\lambda\lambda}dz_j\wedge d\ol z_k+\!
\sum_{\lambda<\mu\le r}d\xi_{\lambda\mu}\wedge d\ol\xi_{\lambda\mu}
-\!\sum_{1\le\mu<\lambda}d\xi_{\mu\lambda}\wedge d\ol\xi_{\mu\lambda},\cr
\Theta_{Q^a}(z_0,\xi_0)&=\sum_{1\le\lambda\le r}a_\lambda
\Theta_{Q_\lambda}(z_0,\xi_0)\cr
&=-\sum_{1\le j,k\le n,\,1\le\lambda\le r}
a_\lambda c_{jk\lambda\lambda}dz_j\wedge d\ol z_k
+\sum_{1\le\lambda<\mu\le r}(a_\lambda-a_\mu)
d\xi_{\lambda\mu}\wedge d\ol\xi_{\lambda\mu}.\cr}
$$
This calculation holds true only at $(z_0,\xi_0)$, but
it shows that we have at every point a decomposition of $\Theta_{Q^a}$ 
in horizontal and vertical parts given by
$$
\leqalignno{
&\qquad\Theta_{Q^a}=\theta_a^H+\theta^V_a,&(\rA.5)\cr
\noalign{\vskip5pt}
&\qquad\theta_a^H(z_0,\xi_0)=-\sum_{j,k,\lambda}
a_\lambda c_{jk\lambda\lambda}dz_j\wedge d\ol z_k
=-\sum_{1\le\lambda\le r}a_\lambda\,
\pi^*\langle\Theta_{T_X,\omega}(e_\lambda),e_\lambda\rangle,&(\rA.6^H)\cr
&\qquad\theta_a^V(z_0,\xi_0)=\sum_{1\le\lambda<\mu\le r}(a_\lambda-a_\mu)
d\xi_{\lambda\mu}\wedge d\ol\xi_{\lambda\mu}.&(\rA.6^V)\cr}
$$
The decomposition is taken here with respect to the $C^\infty$ splitting of 
the exact sequence
$$
0\to T_{Y/X}\to T_Y\to \pi^* T_X\to 0,\qquad Y:=F(E^*)\leqno(\rA.7)
$$
provided by the Chern connection $\nabla$ of $(E,h)\,$; horizontal directions
are those coming from flags associated with $\nabla$-parallel frames.
In order to express $(\rA.6^H)$ in a more intrinsic way at an arbitrary point 
$(z,\xi)\in Y$, we have to replace $(e_\lambda(z))$ by
the orthonormal frame $(\widehat e_\lambda(z,\xi))$ associated with the
flag $\xi\,$; such frames are not unique, actually they are defined up 
to the action of $(S^1)^r$, but such a change does not affect the
expression of $\theta_a^H$. We then get the intrinsic formula
$$
\leqalignno{
\qquad\theta_a^H(z,\xi)&=-\sum_{1\le\lambda\le r}a_\lambda\,
\pi^*\big\langle\Theta_{T_X,\omega}(\widehat e_\lambda(z,\xi)),
\widehat e_\lambda(z,\xi)\big\rangle\cr
&=-\sum_{1\le\lambda\le r}a_\lambda
\sum_{1\le j,k\le n,\,1\le\sigma,\tau\le r}c_{jk\sigma\tau}(z)\,
\widehat e_{\lambda\sigma}(z,\xi)\,
\overline{\widehat e_{\lambda\tau}(z,\xi)}\,dz_j\wedge d\overline z_k
&(\rA.8)\cr}
$$
where we put
$$
\widehat e_\lambda(z,\xi)=\sum_{1\le\sigma\le r}
\widehat e_{\lambda\sigma}(z,\xi)\,e_\sigma(z)
$$
(the coefficients $\widehat e_{\lambda\sigma}(z,\xi)$ can be computed from
(\rA.4)). Moreover, since $\theta_a^V$ and $\Theta_{Q^a}$
have the same restriction to the fibers of $Y\to X$, we conclude that
$\theta_a^V$ is in fact unitary invariant along the fibers (the tautological
metric of $Q^a$ clearly has this property). Let us consider the vertical
and normalized unitary invariant relative volume form $\eta$ of $Y\to X$
given by
$$
\eta(z_0,\xi_0)=\bigwedge_{1\le\lambda<\mu\le r}
i\,d\xi_{\lambda\mu}\wedge d\ol\xi_{\lambda\mu}
\qquad\hbox{at $(z_0,\xi_0)$.}\leqno(\rA.9)
$$
Let $N=r(r-1)/2$ be the fiber dimension. For a strictly dominant
weight $a$, i.e.\ $a_1>a_2>\ldots>a_r$, the line bundle $Q^a$ is
relatively ample with respect to the projection $\pi:Y=F(E^*)\to X$,
and $i\theta_a^V$ induces a K\"ahler form on the fibers. Formula $(\rA.6^V)$
shows that the corresponding volume form is
$$
(i\theta_a^V)^N=N!\prod_{1\le\lambda<\mu\le r}(a_\lambda-a_\mu)~\eta.
$$
\vskip-4pt

\claim A.10. Curvature formulas|Consider as above $Q^a\to Y:=F(E^*)$.
Then$\,:$
\vskip3pt
\item{\rm (a)} The curvature form of $Q^a$ is given
by $\Theta_{Q^a}=\theta_a^H+\theta_a^V$ where the horizontal part is given by
$$
\theta_a^H=-\sum_{1\le\lambda\le r}a_\lambda\,
\pi^*\langle\Theta_{T_X,\omega}(\widehat e_\lambda),\widehat e_\lambda\rangle
$$
and the vertical part by
$$
\theta_a^V(z_0,\xi_0)=\sum_{1\le\lambda<\mu\le r}(a_\lambda-a_\mu)
d\xi_{\lambda\mu}\wedge d\ol\xi_{\lambda\mu}
$$
in normal coordinates at any point $(z_0,\xi_0)$.\vskip3pt
\item{\rm (b)} The relative canonical bundle $K_{Y/X}$ is isomorphic
with $Q^\rho$ for the $($anti-dominant$\,)$ canonical
weight $\rho_\lambda=2\lambda-r-1$, 
$1\le\lambda\le r$. For any positive definite $(1,1)$-form $\omega$ 
on $X$ we have
$$
\eqalign{
i\ddbar\eta\wedge\pi^*\omega^{n-1}
&=-i\theta_\rho^H\wedge \eta\wedge\pi^*\omega^{n-1}\cr
&=\sum_{1\le\lambda\le r}\rho_\lambda\,\pi^*\langle 
i\Theta_{T_X,\omega}(\widehat e_\lambda),\widehat e_\lambda\rangle
\wedge\eta\wedge\pi^*\omega^{n-1}.\cr}
$$
\endclaim

\noindent{\it Proof.} (a) follows entirely from the previous discussion.

\noindent (b) The formula for the canonical weight is a classical result
in the theory of flag varieties. As $(i\theta_a^V)^N$ and $\eta$ are 
proportional for $a$ strictly dominant, we compute instead
$$
\ddbar(\theta_a^V)^N=N\,(\theta_a^V)^{N-1}\wedge\ddbar\theta_a^V+
N(N-1)\,(\theta_a^V)^{N-2}\wedge\partial\theta_a^V\wedge\dbar\theta_a^V,
$$
and for this, we use a Taylor expansion of order 2 at $(z_0,\xi_0)$.
Since $\Theta_{Q^a}$ is closed, we have $\ddbar\theta_a^V=-\ddbar\theta_a^H$,
hence
$$
\ddbar\theta_a^V=\ddbar
\sum_{1\le\lambda\le r}a_\lambda
\sum_{1\le j,k\le n,\,1\le\sigma,\tau\le r}c_{jk\sigma\tau}(z)\,
\widehat e_{\lambda\sigma}(z,\xi)\,
\overline{\widehat e_{\lambda\tau}(z,\xi)}\,dz_j\wedge d\overline z_k,
$$
and we have similar formulas for $\partial(\theta_a^V)$ and $\dbar(\theta_a^V)$.
When taking $\partial$, $\dbar$ and $\ddbar$ we need only consider the
differentials in $\xi$, otherwise we get terms 
$\Lambda^{\ge 3}(dz,d\overline z)$ of degree at least $3$ in the $dz_j$ or 
$d\overline z_k$ and the wedge product of these
with $\pi^*\omega^{n-1}$ is zero. For the same reason, $\partial\theta_a^V
\wedge\dbar\theta_a^V$ will not contribute to the result since it produces 
terms of  degree${}\ge 4$ in $dz_j$, $d\overline z_k$. Formula (\rA.4) gives
$$
\eqalign{
\widehat e_{\lambda\sigma}=\bOne_{\{\lambda=\sigma\}}\Big(
1&-{1\over 2}\sum_{1\le\mu<\lambda}|\xi_{\mu\lambda}|^2
-{1\over 2}\sum_{\lambda<\mu\le r}|\xi_{\lambda\mu}|^2\Big)\cr
&+\bOne_{\{\lambda<\sigma\}}\xi_{\lambda\sigma}
-\bOne_{\{\sigma<\lambda\}}\Big(
\overline\xi_{\sigma\lambda}+\sum_{\mu>\lambda}\xi_{\lambda\mu}
\overline\xi_{\sigma\mu}\Big)+O(|z|^2+|\xi|^3).\cr}
$$
Notice that we do not need to look at the terms $O(|z|)$, $O(|z|^2)$ as 
they will produce no contribution at $(z_0,\xi_0)$. From this we infer
$$
\eqalign{
\widehat e_{\lambda\sigma}&\overline{\widehat e_{\lambda\tau}}
=\bOne_{\{\lambda=\sigma=\tau\}}\Big(
1-\sum_{1\le\mu<\lambda}|\xi_{\mu\lambda}|^2
-\sum_{\lambda<\mu\le r}|\xi_{\lambda\mu}|^2\Big)\cr
&+\bOne_{\{\lambda=\tau<\sigma\}}\xi_{\lambda\sigma}
-\bOne_{\{\sigma<\lambda=\tau\}}\overline\xi_{\sigma\lambda}
+\bOne_{\{\lambda=\sigma<\tau\}}\overline\xi_{\lambda\tau}
-\bOne_{\{\tau<\lambda=\sigma\}}\xi_{\tau\lambda}
\cr
&+\bOne_{\{\sigma,\tau>\lambda\}}\xi_{\lambda\sigma}\overline\xi_{\lambda\tau}
+\bOne_{\{\sigma,\tau<\lambda\}}\xi_{\tau\lambda}\overline\xi_{\sigma\lambda}
+\bOne_{\{\tau<\lambda<\sigma\}}\xi_{\lambda\sigma}\xi_{\tau\lambda}
+\bOne_{\{\sigma<\lambda<\tau\}}\overline\xi_{\sigma\lambda}
\overline\xi_{\lambda\tau}\cr
&-\sum_{1\le\mu\le r}
\bOne_{\{\sigma<\lambda=\tau<\mu\}}\xi_{\lambda\mu}\overline\xi_{\sigma\mu}+
\bOne_{\{\tau<\lambda=\sigma<\mu\}}\xi_{ \tau\mu}\overline\xi_{\lambda\mu}~~
\mod(|z|^2,|\xi|^3).\cr}
$$
In virtue of (\rA.7), only ``diagonal terms'' of the form 
$d\xi_{\lambda\mu}\wedge d\overline\xi_{\lambda\mu}$ in the $\ddbar$ of this
expression can contribute to $(\theta_a^V)^{N-1}\wedge\ddbar\theta_a^V$, 
all others vanish at $z=\xi=0$. The useful terms are thus 
$$\eqalign{
\ddbar&\Big(\sum_{1\le\lambda\le r}a_\lambda\sum_{1\le\sigma,\tau\le r}
c_{jk\sigma\tau}\widehat e_{\lambda\sigma}\overline{\widehat e_{\lambda\tau}}
dz_j\wedge d\overline z_k\Big)=\hbox{(unneeded terms)}{\,}+\cr
&+\sum_{1\le\lambda<\mu\le r}\kern-4pt
\big(\!-a_\mu c_{jk\mu\mu}-a_\lambda c_{jk\lambda\lambda}
+a_\lambda c_{jk\mu\mu}+a_\mu c_{jk\lambda\lambda}\big)
d\xi_{\lambda\mu}\wedge d\overline\xi_{\lambda\mu}
\wedge dz_j\wedge d\overline z_k\cr
&=\sum_{1\le\lambda<\mu\le r}
(a_\lambda-a_\mu)(c_{jk\mu\mu}-c_{jk\lambda\lambda})\,
d\xi_{\lambda\mu}\wedge d\overline\xi_{\lambda\mu}
\wedge dz_j\wedge d\overline z_k+\hbox{(unneeded)}.\cr}
$$
From this we infer
$$
\ddbar(\theta_a^V)^N\wedge\pi^*\omega^{n-1}=
(\theta_a^V)^N\wedge
\sum_{1\le j,k\le n,\,1\le\lambda\le r}(2\lambda-1-r)\,c_{jk\lambda\lambda}\,
dz_j\wedge d\overline z_k\wedge\pi^*\omega^{n-1},
$$
in fact, the coefficient of $c_{jk\lambda\lambda}$ is the number $(\lambda-1)$ 
of indices${}<\lambda$ (coming from the term 
$(a_\lambda-a_\mu)c_{jk\mu\mu}$ above)
minus the number $r-\lambda$ of indices${}>\lambda$ (coming from the term
$-(a_\lambda-a_\mu)c_{jk\lambda\lambda}$). Formula \rA.10~(b) follows.\qed

\claim A.11. Bochner formula|Assume that $X$ is a compact complex manifold 
possessing a balanced metric, i.e.\ a positive smooth $(1,1)$-form
$\omega=i\sum_{1\le j,k\le n}\omega_{jk}(z)\,dz_j\wedge d\overline z_k$ such 
that $d\omega^{n-1}=0$. Assume also that for some dominant weight $a$ 
$(a_1\ge\ldots\ge a_r\ge 0)$, the $\bR$-line bundle $Q^a$ is 
pseudoeffective on~$Y:=F(E^*)$, i.e.\ that there exists
a quasi-plurisubharmonic function $\varphi$ such that
$i(\Theta_{Q^a}+\ddbar\varphi)\ge 0$ on $Y$. Then
$$
\int_Y \big(i\partial\varphi\wedge\dbar\varphi
-i\theta_{a-\rho}^H\big)e^\varphi\,\eta\wedge\pi^*\omega^{n-1}\le 0,
$$
or equivalently
$$
\int_Y \Big(i\partial\varphi\wedge\dbar\varphi
+\sum_{1\le\lambda\le r}(a_\lambda-\rho_\lambda)
\langle i\Theta_{TX,\omega}(\widehat e_\lambda),\widehat e_\lambda\rangle\Big)
e^\varphi\,\eta\wedge\pi^*\omega^{n-1}\le 0.
$$
\endclaim

\noindent{\it Proof.} The idea is to use
the $\ddbar$-formula
$$\eqalign{
\int_Y i\ddbar&(e^\varphi)\wedge \eta\wedge\pi^*\omega^{n-1}-
e^\varphi\wedge i\ddbar\eta\wedge\pi^*\omega^{n-1}\cr
&=\int_Y d\Big(i\dbar(e^\varphi)\wedge \eta\wedge\pi^*\omega^{n-1}+
e^\varphi\, i\partial\eta\wedge\pi^*\omega^{n-1}\Big)=0\cr}
$$
which follows immediately from Stokes. We get
$$
\int_Y (i\ddbar\varphi +i\partial\varphi\wedge\dbar\varphi)e^\varphi\wedge 
\eta\wedge\pi^*\omega^{n-1}-
e^\varphi\, i\ddbar\eta\wedge\pi^*\omega^{n-1}=0.\leqno(\rA.12)
$$
Now, $i\ddbar\varphi\ge -i\Theta_{Q^a}$ in the sense of currents, and 
therefore by \rA.10~(a,b) we obtain
$$
i\ddbar\varphi\wedge\eta\wedge\pi^*\omega^{n-1}
-i\ddbar\eta\wedge\pi^*\omega^{n-1}
\ge (-i\theta_a^H+i\theta_\rho^H)\wedge\eta\wedge\pi^*\omega^{n-1}.\leqno(\rA.13)
$$
The combination of (\rA.12) and (\rA.13) yields the inequality of 
Corollary~\rA.11.\qed
\medskip

The parallel transport operators of $(E,h)$ can be considered to operate 
on the global flag variety $Y=F(E^*)$ as follows. For any piecewise smooth 
path $\gamma:[0,1]\to X$, 
we get a (unitary) hermitian isomorphism $\tau_\gamma:E_p\to E_q$ where
$p=\gamma(0)$, $q=\gamma(1)$. Therefore $\tau_\gamma$ induces an
isomorphism $\widetilde \tau_\gamma:F(E^*_p)\to F(E^*_q)$
of the corresponding flag varieties, and an isomorphism over
$\widetilde\tau_\gamma$ of the tautological line bundles~$Q^a$. Given
a local $C^\infty$ vector field $v$ on an open open set $U\subset X$,
there is a unique horizontal lifting $\widetilde v$ of $v$ to a
$C^\infty$ vector field on $\pi^{-1}(U)\subset Y$, where
horizontality refers again to~$\nabla=\nabla_{E,h}$. Now, the flow of
$\widetilde v$ consists of parallel transport operators along the
trajectories of~$v$. By definition, $h$ is invariant by
parallel transport, therefore the associated hermitian metric $h_a$
on each line bundle $Q^a$ is also invariant. Another metric 
$h_{a,\varphi}=h_ae^{-\varphi}$ is invariant if and only if the weight
function $\varphi$ is invariant by the flows of all such vector fields
$\widetilde v$ on $Y$, that is if $d\varphi(\zeta)=0$ for all
horizontal vector fields $\zeta\in T_Y$.

\claim A.14. Theorem|Let $E\to X$ be a holomorphic vector bundle of rank $r$
over a compact complex manifold~$X$. Assume that $X$ is equipped with a 
hermitian metric $\omega$ and $E$ with a hermitian structure $h$ 
such that $B:=\Tr_\omega(i\Theta_{E,h})\ge 0$. At each point $z\in X$, let
$$
0\le b_1(z)\le \ldots\le b_r(z)
$$
be the eigenvalues of $B(z)$ with respect to~$h(z)$. Finally,
let $Q^a$ be a pseudoeffective $\bR$-line bundle on $Y:=F(E^*)$ associated 
with a dominant weight $a_1\ge\ldots\ge a_r\ge 0$, and let $\varphi$ be a 
quasi-plurisubharmonic function on $Y$ such that 
$i(\Theta_{Q^a}+\ddbar\varphi)\ge 0$. Then
\vskip3pt
\item{\rm (a)} The function $\psi(z)=\sup_{\xi\in F(E^*_z)}\varphi(z,\xi)$
is constant and $b_\lambda\equiv 0$ as soon as $a_\lambda>0$, and in particular
$B\equiv 0$ if $a_r>0$.
\vskip3pt
\item{\rm (b)} Assume that $B\equiv 0$. Then the function $\varphi$ must 
be invariant by parallel transport on~$Y$.
\vskip3pt
\endclaim

\noindent{\it Proof.} Since our hypotheses are invariant by a conformal 
change on the metric $\omega$, we can assume by Gauduchon [Gau77] that
$\ddbar\omega^{n-1}=0$.
\medskip

\noindent
(a) Notice that if $a$ is integral and $\varphi$ is 
given by a holomorphic section of $Q^a$, then $e^\varphi$ is the
square of the norm of that section with respect to~$h$, and $e^\psi$ is
the sup of that norm on the fibers of $Y\to X$. In general, formula A.10~(a)
shows that
$$
i\ddbar\varphi(z,\xi)\ge -i\theta_a^H(z,\xi)-i\theta_a^V(z,\xi),
$$
hence
$$
i\ddbar^H\!\!\varphi(z,\xi)\wedge\omega^{n-1}(z)\ge \sum_{1\le\lambda\le r}a_\lambda\,
\langle i\Theta_{TX,\omega}(\widehat e_\lambda),\widehat e_\lambda\rangle(z,\xi)
\wedge \omega^{n-1}(z)\leqno(\rA.15)
$$
where $i\ddbar^H\!\!\varphi$ means the restriction of
$i\ddbar\varphi$ to the horizontal directions in $T_Y$. By taking the
supremum in $\xi$, we conclude from standard arguments of subharmonic function
theory that
$$
\Delta_\omega\psi(z)\ge \sum_{1\le\lambda\le r}a_\lambda b_\lambda(z),
$$
since the right hand side is the minimum of the coefficient of the
$(n,n)$-form occurring in the RHS of (\rA.15). Therefore $\psi$ is 
$\omega$-subharmonic and so must be constant on $X$ by Aronszajn 
[Aro57]. It follows that $b_\lambda\equiv 0$ whenever $a_\lambda>0$,
in particular $B\equiv 0$ if $a_r>0$.
\medskip

\noindent(b) Under the assumption $B=\Tr_\omega\Theta_{E,h}\equiv 0$, the 
calculations made in the course of the proof of A.10~(b) imply that
$$
\partial\eta\wedge\pi^*\omega^{n-1}=0,\qquad 
\ddbar\eta\wedge\pi^*\omega^{n-1}=0.
$$
By the proof of the Bochner formula A.11 (the fact that 
$\ddbar\omega^{n-1}=0$ is enough here), we get
$$
0\le \int_Y i\partial\varphi\wedge\dbar\varphi\wedge \eta\wedge\pi^*
\omega^{n-1}\le 0,
$$
and we conclude from this that the horizontal derivatives
$\partial^H\!\varphi$ vanish. Therefore $\varphi$ is invariant by parallel
transport.\qed
\medskip

In the vein of Criterion~1.1, we have the following additional statement.

\claim A.16. Proposition|Let $X$ be a compact K\"ahler manifold. Then
$X$ is projective and rationally connected if and only if none of the
$\bR$-line bundles $Q^a$ over $Y=F(T_X^*)$ is pseudoeffective
for weights $a\ne 0$ with $a_1\ge\ldots\ge a_r\ge 0$.
\endclaim

\noindent{\it Proof.} If $X$ is projective rationally connected
and some $Q^a$, $a\ne 0$, is pseudoeffective, we obtain a
contradiction with Theorem A.14 by pulling-back $T_X$ and $Q^a$ via a
map $f:\bP^1\to X$ such that $E=f^*T_X$ is ample on $\bP^1$ (as $B>0$
in this circumstance).

Conversely, if the $\bR$-line bundles $Q^a$, $a\ne 0$, are not 
pseudoeffective on $Y=F(T^*_X)$, we obtain
by taking $a_1=\ldots=a_p=1$, $a_{p+1}=\ldots=a_n=0$ that
$\pi_*Q^a=\Omega^p_X$. Therefore $H^0(X,\Omega^p_X)=0$ and all 
invertible subsheaves $\cF\subset\Omega^p_X$ are not 
pseudoeffective for $p\ge 1$. Hence $X$ is projective (take $p=2$ and
apply Kodaira [Kod54]) and rationally connected by Criterion 1.1~(b).\qed
\vskip20pt

\section References|

\bibitem[Aro57]&Aronszajn, N.:& A unique continuation theorem for
solutions of elliptic partial differential equations or inequalities
of second order.& J.~Math.\ Pures Appl.\ {\bf 36}, 235--249 (1957)&

\bibitem[Aub76]&Aubin, T.:& Equations du type Monge-Amp\`ere sur les
vari\'et\'es k\"ahle\-riennes compactes.& C.\ R.\ Acad.\ Sci.\ Paris
Ser.\ A {\bf 283}, 119--121 (1976)$\,$; Bull.\ Sci.\ Math.\
{\bf 102}, 63--95 (1978)&

\bibitem[Bea83]&Beauville, A.:& Vari\'et\'es k\"ahleriennes dont la
premi\`ere classe de Chern est nulle.&
J.\ Diff.\ Geom.\ {\bf 18}, 775--782 (1983)& 

\bibitem[Ber55]&Berger, M.:& Sur les groupes d'holonomie des vari\'et\'es
\`a connexion affine des vari\'et\'es riemanniennes.& Bull.\ Soc.\ Math.\
France {\bf 83}, 279--330 (1955)& 

\bibitem[Bes87]&Besse, A.L.:& Einstein manifolds.& Springer-Verlag,
New York, 1987&
  
\bibitem[Bis63]&Bishop, R.:& A relation between volume, mean curvature and
diameter.& Amer.\ Math.\ Soc.\ Not.\ {\bf 10}, p.$\,$364 (1963)& 

\bibitem [BY53]&Bochner, S., Yano, K.:& Curvature and Betti numbers.& 
Annals of Mathematics Studies, No.~32, Princeton University Press, 
Princeton, N.~J., ix$+190\,$pp (1953)&

\bibitem[Bog74a]&Bogomolov, F.A.:& On the decomposition of K\"ahler
manifolds with trivial canonical class.& Math.\ USSR Sbornik {\bf 22},
580--583 (1974)& 

\bibitem[Bog74b]&Bogomolov, F.A.:& K\"ahler manifolds with trivial canonical
class.& Izves\-tija Akad.\ Nauk {\bf 38}, 11--21 (1974)& 

\bibitem[BDPP]&Boucksom, S., Demailly, J.-P., Paun, M., Peternell, T.:&
The pseudo-effective cone of a compact K\"ahler manifold and varieties 
of negative Kodaira dimension.& arXiv: 0405285 [math.AG]$\,$;
J.\ Alg.\ Geometry {\bf 22}, 201--248 (2013)& 

\bibitem[Bru10]&Brunella, M.:& On K\"ahler surfaces with semipositive
Ricci curvature.& Riv.\ Math.\ Univ.\ Parma (N.S.), {\bf 1}, 441--450 (2010)&

\bibitem[Cam92]&Campana, F.:& Connexit\'e rationnelle des vari\'et\'es 
de Fano.& Ann.\ Sci.\ Ec.\ Norm. \ Sup.\ {\bf 25}, 539--545 (1992)&

\bibitem[Cam95]&Campana, F.:& Fundamental group and positivity of cotangent
bundles of compact K\"ahler manifolds.& J.\ Alg.\ Geom.\ {\bf 4},
487--502 (1995)&

\bibitem[CPZ03]&Campana, F., Peternell, Th., Zhang, Qi:&
On the Albanese maps of compact K\"ahler manifolds.&
Proc.\ Amer.\ Math.\ Soc.\ {\bf 131}, 549--553 (2003)&

\bibitem[CG71]&Cheeger, J., Gromoll, D.:& The splitting theorem for
manifolds of nonnegative Ricci curvature.& J.\ Diff.\ Geom.\ {\bf 6},
119--128 (1971)& 

\bibitem[CG72]&Cheeger, J., Gromoll, D.:& On the structure of complete
manifolds of nonnegative curvature.& Ann.\ Math.\ {\bf 96}, 413--443
(1972)&

\bibitem[DPS93]&Demailly, J.-P., Peternell, T., Schneider, M.:& K\"ahler
manifolds with numerically effective Ricci class.& Compositio Math.\ 
{\bf 89}, 217--240 (1993)& 

\bibitem[DPS94]&Demailly, J.-P., Peternell, T., Schneider, M.:& Compact
complex manifolds with numerically effective tangent bundles.& J.\ Alg.\
Geom.\ {\bf 3}, 295--345 (1994)& 

\bibitem[DPS96]&Demailly, J.-P., Peternell, T., Schneider, M.:& Compact
K\"ahler manifolds with hermitian semipositive anticanonical
bundle.& Compositio Math.\ {\bf 101}, 217--224 (1996)&

\bibitem[DR52]&de Rham, G.:& Sur la reductibilit\'e d'un espace de Riemann.& 
Comment.\ Math.\ Helv.\ {\bf 26} 328--344 (1952)& 

\bibitem[Gau77]&Gauduchon, P.:& Le th\'eor\`eme de l'excentricit\'e nulle.&
C.~R.\ Acad.\ Sci.\ Paris {\bf 285}, 387--390 (1977)&

\bibitem [GHS01]&Graber, T., Harris, J.,Starr, J.:& Families of
rationally connected varieties.& Journal of the Amer.\ Math.\ Soc.\
{\bf 16}, 57--67 (2003)&

\bibitem [K\kern-0.6pt M\kern-1.6pt M92]&Koll\'ar, J., Miyaoka,Y., 
Mori,S:& Rationally connected varieties.& J. \ Alg.\ {\bf 1}, 429--448 (1992)&

\bibitem[Kob81]&Kobayashi, S.:& Recent results in complex differential
geometry.& Jber.\ dt.\ Math.-Verein.\ {\bf 83}, 147--158 (1981)&

\bibitem[Kob83]&Kobayashi, S.:& Topics in complex differential geometry.&
In DMV Seminar, Vol.~3., Birkh\"auser 1983& 

\bibitem[Kod54]&Kodaira, K.:& On K\"ahler varieties of restricted type.&
Ann.\ of Math.\ {\bf 60}, 28--48 (1954)&

\bibitem[Kol96]&Koll\'ar, J.:& Rational Curves on Algebraic Varieties.& 
Ergebnisse der Mathematik und ihrer Grenzgebiete, 3.~Folge, Band 32,
Springer, 1996&

\bibitem[Lic67]&Lichnerowicz, A.:& Vari\'et\'es k\"ahleriennes et
premi\`{e}re classe de Chern.&
J.\ Diff.\ Geom.\ {\bf 1}, 195--224 (1967)& 

\bibitem[Lic71]&Lichnerowicz, A.:& Vari\'et\'es K\"ahl\'eriennes \`a
premi\`ere classe de Chern non n\'egative et vari\'et\'es riemanniennes \`a
courbure de Ricci g\'en\'eralis\'ee non n\'egative.& J.\ Diff.\ Geom.\ {\bf
6}, 47--94 (1971)& 

\bibitem[Pau12]&P\u{a}un, M.:& Relative adjoint transcendental classes and 
the Albanese map of compact K\"ahler manifolds with nef Ricci classes.&
arXiv: 1209.2195 &

\bibitem[Pet06]&Peternell, Th.:& Kodaira dimension of subvarieties II.& Intl.\ 
J.\ Math.\ {\bf 17}, 619--631 (2006)&

\bibitem[PS98]&Peternell, Th., Serrano, F.:& Threefolds with anti canonical bundles& Coll.\ 
Math.\ {\bf 49}, 465--517 (1998)&

\bibitem [Ued82]&Ueda, T.:& On the neighborhood of a compact complex curve 
with topologically trivial normal bundle.& J.\ Math.\ Kyoto Univ.\
{\bf 22}, 583--607 (1982/83)&

\bibitem[Yau78]&Yau, S.T.:& On the Ricci curvature of a complex K\"ahler
manifold and the complex Monge-Amp\`ere equation~I.& Comm.\ Pure and Appl.\
Math.\ {\bf 31}, 339--411 (1978)&

\bibitem[Zha96]&Zhang, Qi:& On projective manifolds with nef anticanonical
bundles.& J.~Reine Angew.\ Math.\ {\bf 478}, 57--60 (1996)&

\bibitem[Zha05]&Zhang, Qi:& On projective varieties with nef anticanonical
divisors.& Math.\ Ann.\ {\bf 332}, 697--703 (2005)&
\vskip15pt

{\parindent=0cm
(version of October 1, 2012, revised on July 16, 2013,
minor  additional corrections made on February 3, 2018,
printed on \today) \vskip15pt
\bigskip
{\bf Fr\'ed\'eric Campana}\hfill\break
Universit\'e de Lorraine, Institut \'Elie Cartan,\hfill\break
UMR 7502 du CNRS, BP~70239\hfill\break
54506 Vand{\oe}uvre-l\`es-Nancy Cedex, France\hfill\break
{\it e-mail\/}: frederic.campana@univ-lorraine.fr\vskip6pt 
{\bf Jean-Pierre Demailly}\hfill\break
Universit\'e de Grenoble I, Institut Fourier,\hfill\break
UMR 5582 du CNRS, BP 74\hfill\break
38402 Saint-Martin d'H\`eres, France\hfill\break
{\it e-mail\/}: demailly@fourier.ujf-grenoble.fr\vskip6pt
{\bf Thomas Peternell}\hfill\break
Universit\"at Bayreuth, Mathematisches Institut\hfill\break D-95440
Bayreuth, Deutschland\hfill\break
{\it e-mail\/}: thomas.peternell@uni-bayreuth.de}


\bye